\documentclass{amsart}

\usepackage{amsmath}
\usepackage{amsfonts}
\usepackage{amsthm}

\usepackage{graphicx} 
\usepackage{amssymb,mathrsfs}
\usepackage{url}

\usepackage{color}

\newtheorem{theorem}{Theorem}[section]
\newtheorem{lemma}[theorem]{Lemma}
\newtheorem{corollary}[theorem]{Corollary}
\newtheorem{proposition}[theorem]{Proposition}

\theoremstyle{definition}
\newtheorem{definition}[theorem]{Definition}
\newtheorem{algorithm}[theorem]{Algorithm}

\theoremstyle{remark}

\numberwithin{equation}{section}

\newcommand{\abs}[1]{\lvert#1\rvert}

\newcommand{\cS}{\mathcal{S}}
\newcommand{\cN}{\mathcal{N}}

\newcommand{\cU}{\mathcal{U}}

\newcommand{\bN}{\mathbb{bN}}
\newcommand{\bR}{\mathbb{R}}
\newcommand{\bQ}{\mathbb{Q}}
\newcommand{\bC}{\mathbb{C}}
\newcommand{\bZ}{\mathbb{Z}}

\newcommand{\sC}{\mathscr{C}}
\newcommand{\sH}{\mathscr{H}}
\newcommand{\sB}{\mathscr{B}}

\newcommand{\CC}{\mathbb{C}}
\newcommand{\ep}{\varepsilon}

\newcommand{\cD}{\mathcal{D}}
\newcommand{\cB}{\mathcal{B}}
\newcommand{\cR}{\mathcal{R}}

\renewcommand{\Re}{\text{Re}}
\renewcommand{\Im}{\text{Im}}

\newcommand{\norm}[1]{\left\| #1\right\|}
\newcommand{\onorm}[1]{\abs{#1}}
\newcommand{\snorm}[1]{\norm{#1}}

\begin{document}

\title[Poly-time computability of hyperbolic Julia sets]
{A dynamical algorithm to compute hyperbolic Julia sets in polynomial time}
 
\author[S. Boyd]{Suzanne Boyd}
\address{Department of Mathematical Sciences\\
University of Wisconsin Milwaukee\\
PO Box 413\\
Milwaukee, WI 53201, 
USA}
\email{sboyd@uwm.edu, ORCID: 0000-0002-9480-4848}

\author[C. Wolf]{Christian Wolf}
\address{Department of Mathematics and Statistics\\
Mississippi State University\\
Starkville, MS 39759, USA}
\email{cwolf@math.msstate.edu, ORCID: 0000-0002-7976-3574.}
\thanks{C.W. was partially supported by grants from the Simons Foundation (SF-MPS-CGM-637594 and SFI-MPS-TSM-00013897) and PSC-CUNY (TRADB-55-67482, TRADB-56-68594).}

\date{\today}


\begin{abstract}
    Hyperbolic Julia sets of  complex polynomials are known to be computable in polynomial time due to pioneering work of Braverman in 2005 \cite{Braverman2005}. In this paper, we present an alternative method for establishing poly-time computability of hyperbolic Julia sets, which allows us to establish, via a new algorithm, lower computability of the hyperbolicity locus of polynomials  of a fixed degree. We first adapt our recently developed algorithms for the computability of polynomial skew products \cite{BoydWolf-Skew1}, and then apply a refinement that allows us to establish poly-time computation of hyperbolic Julia sets. Finally, we derive lower computability of the hyperbolicity locus via an adapted lattice/refinement search algorithm. In contrast to \cite{Braverman2005}, our approach is dynamical in nature and does not rely on techniques unique to complex analysis. 
\end{abstract}

\maketitle


\section{Introduction}
\label{sec:Intro}

Analyzing the computability and non-computability of dynamically-defined sets and invariants has been an important area of research during the last 25 years, see, e.g., \cite{BBRY-2011,Binder-et-al-2024,BY2009,Burr-Wolf-2024,DY2018}. This area is often referred to as ``Computability in Dynamical Systems". Particular attention has been given to the computability of Julia sets in one-dimensional complex dynamics, see, e.g., \cite{BY2009,Burr-Wolf-2024} and the references therein. One of the first striking results in this area was due to Braverman \cite{Braverman2005} who showed that hyperbolic Julia sets are computable in polynomial time; also called poly-time computable.  Braverman's proof heavily relies  on complex-analytic tools, such as working with the Poincare metric and applying Pick's theorem. In this paper, we develop an alternative algorithm for the poly-time computability of hyperbolic Julia sets. The main difference from Braverman's approach is that our algorithm relies on dynamical methods. One motivation of this paper is to develop a dynamical computability framework in the hope that it can be generalized to answer computability questions for dynamical systems beyond the complex-analytic case. 

Our methods in part rely on ideas we recently developed to study computability for polynomial skew products of $\bC^2$ \cite{BoydWolf-Skew1}. Besides reproving Braverman's result, we also obtain the lower computability of the hyperbolicity locus of polynomials with a fixed degree $d\geq 2$, which follows analogously to our proof for polynomial skew products \cite{BoydWolf-Skew1}.

Braverman \cite{Braverman2005} established polynomial-time computability of a hyperbolic  Julia set in time $O(N\cdot M(N)),$ where $M(N)$ is the time complexity of multiplying two $N$-bit numbers, with best-known estimate $M(N)=O(N\log N \log \log N)$. Our main theorem is the same, though based on an approach that is more accessible to generalization. 

\begin{theorem} \label{thm:J-computable}
    The Julia set of a hyperbolic polynomial of degree $d\geq 2$ is computable to precision $2^N$ in time $O(N\cdot M(N))$. 
\end{theorem}
Roughly speaking, our algorithm can be summarized by the following two steps. \underline{Step 1}: We apply the techniques of \cite{BoydWolf-Skew1} to compute a covering of $J$ by finitely many sufficiently small ``reference" boxes, so that every reference box has a non-empty intersection with $J$. This succeeds if $p$ is hyperbolic. We note that Step 1 is solely based on dynamical properties including the invariance of the Julia set and its characterization as the expanding chain recurrent component.  Moreover, for any given hyperbolic polynomial $p$, Step 1 needs to be performed only once and thus does not contribute to the computational complexity of the precision of the computation of the Julia set of $p$. \underline{Step 2}: If Step 1 is successful (which it is if $p$ is hyperbolic), we consider a small test box $B$ of side-length $2^{-N}$, with center $z$ where $z$ is contained in (or sufficiently close to) one of the reference boxes. The goal is to decide in polynomial (in fact, linear) complexity of $N$ whether or not the test box intersects the Julia set. To accomplish this we compute underestimates, respectively overestimates,
of $p^k(B)$ in terms of balls $B(p^k(z),c_{u/o} (p^k)'(z))$ for some constants $c_u<c_o$ which are independent of $N$.
At each iterate we check if the underestimate of $p^k(B)$ contains a reference box, in which 
case the algorithm reports $1$, that is, $J$ intersects the test box $B$, or if the overestimate of $p^k(B)$ is disjoint from all reference boxes, in which case the algorithm reports $0$, that is, $B$ does not intersect $J$. We prove that this algorithm terminates if $p$ is hyperbolic, and that the number of required iterates is bounded by a constant times $N$. One of the ingredients to make Step 2 work is a distortion estimate that is used to control the diameter of the iterates of $B$  in terms of the derivative of $p^k$ at $z$. While we work here for simplicity with Koebe’s distortion theorem, it should be noted that similar distortion results are known to exist for many smooth uniformly hyperbolic systems.

Because Step 1 of our algorithm requires establishing a neighborhood on which the map is expanding, it implies a corollary for the hyperbolicity locus of polynomials. More precisely,
for $d\geq 2$  let $\mathscr{A}_d$ be the subset of hyperbolic polynomials in 
$\{ z^d + a_{d-2}z^{d-2} + \cdots + a_1 z + a_0: a_\ell \in \CC \} \cong \CC^{d-1}$ parametrized by the coefficients of the polynomial. We have (just as we showed for polynomial skew products, \cite{BoydWolf-Skew1}):
\begin{corollary} \label{cor:semi-decidable}
    Hyperbolicity of polynomials of fixed degree $d$ is a semi-decidable property on $\CC^{d-1}$. More precisely, there exists a Turing machine that, on input of an oracle of $c\in \CC^{d-1}$, halts if $p_c$ is hyperbolic, and runs forever if $p_c$ is not hyperbolic. 
\end{corollary}

Finally, this corollary, and the Algorithm 4.1 and Theorem 1.3 of \cite{BoydWolf-Skew1} restricted to one-dimensional polynomials, yields the following.

\begin{theorem} \label{thm:hyp-locus}
    Let $d\geq 2$. The  locus $\mathscr{A}_d$ of hyperbolic polynomials of degree $d$ is lower semi-computable. 
\end{theorem}
In the above result, lower semi-computability roughly-speaking means that there is an algorithm which outputs a sequence of balls whose union coincides with $\mathscr{A}_d$. This sequence is, in general, infinite, and does not provide  at any given time information about the size of the missing points in $\mathscr{A}_d$.

This paper is organized as follows. In Section~\ref{sec:prelim} we provide the necessary background from computability theory and complex dynamics. We also derive some consequences of the general Koebe distortion theorem that we need for our main theorem.
Since we consider holomorphic maps,  we use Koebe's distortion theorem for simplicity. In fact, since uniformly hyperbolic dynamical systems exhibit bounded distortion, our approach could be adapted in a more general setting by using different results, in the case of more general hyperbolic systems. 

In Section~\ref{sec:results}, we provide our main algorithm and the proofs of the main theorems, in Section~\ref{sec:firststep} describing our ``overhead'' first step of the algorithm, and in Section~\ref{sec:mainresult} detailing the algorithm for increasing precision with $N$, and the proof of Theorem~\ref{thm:J-computable}.

\section{Preliminaries and background material}
\label{sec:prelim}
In this section, we provide an overview of the basic definitions and properties of the material required in this paper.

\subsection{Computability}
\label{sec:compute:basic}
Computability theory provides a tool to describe when mathematical objects can be approximated to any desired accuracy. 
We give key computability results below, but we refer the reader for a more detailed discussion of computability theory to \cite{BoydWolf-Skew1} or in general 
\cite{BBRY-2011, Binder-et-al-2024, Braverman2005, BY2009, BSW2020, Burr-Wolf-2024, GHR2011}. In the following, we use a bit-based computation model (which simply means that information is stored as binary digits), such as a Turing machine (a computer program for our purposes).  One can think of the set of Turing machines as a particular, countable set of functions; we denote $T(x)$ as the output of the Turing machine $T$ based on input $x$.
We start with the computability of points in Euclidean spaces.

\begin{definition}
Let $\ell\in \bN$ and $x\in \bR^\ell$. An \emph{oracle} of $x$ is a function $\phi:\bN\to \bQ^\ell$ such that $\Vert \phi(n)-x\Vert < 2^{-n}$. Moreover, we say $x$ is {\em computable}
if there is a Turing Machine $T=T(n)$ 
which is an oracle of $x$.
\end{definition} 

It is easy to see that rational numbers, algebraic numbers, and some transcendental numbers such as {\it e} and $\pi$ are computable real numbers. However, since the collection of Turing machines is countable, most points in $\bR^\ell$ are not computable.
Identifying $\bC^\ell$ with $\bR^{2\ell}$, the notion of computable points naturally extends to $\bC^\ell$. Next, we define computable functions on Euclidean spaces.

\begin{definition}\label{defcompfunc} 
Let $D\subset \bR^\ell$. A function $f:D\rightarrow\bR^k$ is {\em computable} if there is a Turing machine $T$ so that for any $x\in D$, any oracle $\phi$ for $x$ and any $n\in \bN, \,\, T(\phi,n)$
is a point in $\bQ^k$ so that $\Vert T(\phi,n)-f(x)\Vert<2^{-n}$.
\end{definition}

One of the inputs of the Turing machine $T$ in this definition is an oracle. Specifically, while the Turing machine $T$ in principle has access to an infinite amount of data, it must be able to decide when the approximation $\phi(m)$ of $x$  is sufficiently accurate to perform the computation to precision $2^{-n}$.
We additionally note that 
the input points $x$ of a computable function are not required to be computable.  Next, we extend the notion of computable points to more general spaces, called computable metric spaces.

\begin{definition}\label{def:computable}
Let $(X,d_X)$ be a separable metric space with metric $d_X$, and let $\cS_X=\{s_i: i\in \bN\}\subset X$ be a countable dense subset. We say $(X,d_X,\cS_X)$ is a \textit{computable metric space} if the distance function $d_X(\cdot,\cdot)$ is uniformly computable on $\cS_X\times \cS_X$, that is, if there exists a Turing machine $T=T(i,j,n)$, which on input $i,j,n\in \bN$ outputs a rational number such that
$|d_X(s_i,s_j)-T(i,j,n)|<2^{-n}$.
\end{definition}

The points in $\cS_X$ in this definition are called the \textit{ideal points} of $X$, and $\cS_X$ is the ideal set of the computable metric space. The ideal points take on the role of $\bQ^k$ in $\bR^k$. Clearly, Euclidean spaces are computable metric spaces with $\cS_{\mathbb{R}^k}=\mathbb{Q}^k$.
In a computable metric space $(X,d_X,\cS_X)$,
we say that a ball $B(x,r)$ is an \textit{ideal ball} if $x\in \cS_X$ and $r=2^{-i}$ for some $i \in \bZ$.

 We use the next definition
to study the locus of hyperbolic polynomials of a fixed degree. 
\begin{definition}
\label{defn:LowerComputableSet} Let $(X,d_X,\cS_X)$ be a computable metric space.
Let $D\subset X$ be open. We say $D$ is {\em lower semi-computable} if there exists a Turing machine producing $T =\{(x_i,r_i)\}_{i\in\bN}$ such that $B(x_i,r_i)$ is an ideal ball and 
$D = \cup_{i=1}^\infty B(x_i,r_i)$. 
\end{definition}
We note in the definition of lower semi-computability we do not require an error estimate for how close any finite union of dyadic balls is to the set $D$, just that it converges in the limit. In fact, this does not even require the set $D$ to be bounded. Computability of a set is a much stronger condition than lower semi-computability.
\begin{definition}
Let $(X,d_X,\cS_X)$ be a computable metric space.
An {\em oracle} for $x\in X$ is a function $\phi$ such that on input $n\in \bN$, the output $\phi(n)$ is a natural number so that $d_X(x,s_{\phi(n)})<2^{-n}$.  Moreover, we say $x$ is {\em computable} if there is a Turing machine $T=T(n)$ which is an oracle for $x$.
\end{definition}

We next extend the notion of computable functions on $\bR^k$ to functions between computable metric spaces. 
\begin{definition}\label{def:computablefunction}
Let $(X,d_X,\cS_X)$ and $(Y,d_Y,\cS_Y)$ be computable metric spaces and  $\cS_Y=\{t_i: i\in \bN\}$.  Let $D\subset X$.  A function $f:D\rightarrow Y$ is {\em computable} if there is a Turing machine $T$ such that for any  $x\in D$ and any oracle $\phi$ of $x$, the output $T(\phi,n)$ is a natural number satisfying $d_Y(t_{T(\phi,n)},f(x))< 2^{-n}$.
\end{definition}

We now turn to the computability of compact subsets of computable metric spaces.
Recall the Hausdorff distance between two {\em compact} subsets $A$ and $B$ of a metric space $X$ is the largest distance of a point in one set to the other set:
$$
d_H(A,B)=\max\left\{\max_{a\in A}d(a,B),\max_{b\in B}d(b,A)\right\},
$$
where $d(x,C)=\min\{d_X(x,y): y\in C\}$.  
We write $\cN(S,\delta)$ for the open $\delta$-neighborhood  in the Hausdorff metric about a set $S$, and $B(s,\delta)$ for the open $\delta$-ball about a point $s$.

 For 
    $\sC= \{C\subset X\,\, {\rm compact}\}$ let $\cS_\sC = \cS_\sC(X)$ denote the collection of finite unions of closed ideal balls.
Then (see e.g.\ \cite{BY2009})
$(\sC,d_H,S_\sC)$ is a computable metric space, when $(X,d_X,\cS_X)$ is.

\subsection{Computability in the complex plane}

By default, rather than working with the standard Euclidean norm on $\bC$, we work with the $L^\infty$ norm for which balls become boxes. More precisely, when we write $\norm{\cdot}$ we mean the $L^{\infty}$ 
norm on~$\mathbb{R}^{2}$, so that for a vector $z\in \mathbb{C},$ we have $\snorm{z}=\max\{ \abs{\Re(z)},\abs{\Im(z)}\}. $
Hence, 
$$d(z,w) = \max \{ | \Re(z)-\Re(w)|, |\Im(z)-\Im(w)|  \}, $$ 
and if $z$ and $w$ are in a (closed) box of side-length $r$, then $d(z,w) \leq r$. 

Consider the computable metric space $(\bC,d,\cS)$, and let $C\subset \bC$ be compact.

We provide the definitions needed for a set to be computable, similar to \cite{Braverman-Yampolsky-2008} except underlying these we use the $L^\infty$ metric.

\begin{definition}
\label{defn:2^n-approx}
    A set $\Psi_n$ is a \textit{$2^n$-approximation of $C$} if  $d_H(C,\Psi_n)\leq 2^{-n},$ in the Hausdorff metric.
\end{definition}

\begin{definition} \label{defn:pic-according-to-h_C}
Let $\cS_{n+2}\subset \cD^4$ be the elements whose component coordinate denominators are $2^{-n-2}$, and $\cS= \cup_{n\in\bN} \cS_n$ the set of all ideal points. 
Consider functions $h_C: \{ (n,x): n\in \bN, x\in \cS_{n+2} \} \to \{0,1\}$
satisfying
\begin{equation}\label{eqn:complexsetfunc}
h_C(n,x)=\left\{
\begin{array}{ll}
      1 & d(x,C)< 2^{-n-2} \\
      0 & d(x,C)\geq 2^{-n-1} \\
      0\, {\rm \ or \ }\, 1 & {\rm otherwise: \ } 2^{-n-2} \leq d(x,C) < 2^{-n-1}. \\
\end{array} 
\right.
\end{equation}
Note the domain $\cS_{n+2}$ is a strict subset of the possible ideal points, and since $C$ is compact, $h_C(n,\cdot)$ has finite support: $\sH_n =  \{x\in \cS_{n+2}: h_C(n,x)=1\}$.

The \textit{$L^\infty$-picture drawn according to $h_C(n,x)$} is the set of all squared pixels of radii $2^{-n-3}$ with centers the points in $\sH_n$.
\end{definition}

As shown in \cite{Braverman2005, DY2018}, for each $n$, 
the picture drawn according to $h_C(n,x)$ (which they describe using round pixels) is a $2^{-n}$-approximation of $C$. In fact, a straight-forward argument shows that the same holds true in our setting using square pixels.

\begin{lemma}
\label{lem:picture_precision}
    The $L^\infty$-picture drawn according to $h_C(n,x)$ is a $2^{-n}$ approximation of $C$. Hence  $C$ is computable if a computable $h_C(n,x)$ exists for $C$, $\forall n\in \bN$. 
\end{lemma}

The existence of the computable function $h_C(n,x)$ also implies that $C$ is {computable} under the definition of computable set given in \cite{BoydWolf-Skew1}:  There is a Turing machine $T$ such that for each $n\in\bN$, $T(n)$ produces a set $\Psi_n$ (defined as a finite collection of ideal centers and radii of balls) which is a $2^n$-approximation of $C$.  
It is well-known that these two alternate definitions of computable set are equivalent. The advantage of the function approach is that it facilitates calculation of the algorithm's running time.

\smallskip

\subsection{Poly-time computability.}
So far, we have considered the question  if a point, function, etc.\ is {\it in principle} computable. In the affirmative case it is natural to ask if the corresponding computation can be performed in a ``reasonable" amount of time. This question is addressed in complexity theory; in particular, only computations that can be performed in polynomial time (often abbreviated as poly-time computability) can be performed in practice. A compact set $C$ is \textit{computable in polynomial time} if there exists a computable function $h_C$ satisfying \eqref{eqn:complexsetfunc} and the number of computations required to output the value of $h_C$ is uniformly bounded by a polynomial in $n$.

\begin{definition}\label{def:polyntimefunction}
Let $(X,d_X,\cS_X)$ and $(Y,d_Y,\cS_Y)$ be computable metric spaces and  $\cS_Y=\{t_i: i\in \bN\}$.  Let $D\subset X$, and  $f:D\rightarrow Y$ be {computable} given by Turing machine $T=T(\phi,n)$. 
We use the standard convention for calculating  the running time: querying the oracle $\phi$ at precision $2^{-m}$ takes $m$ units of time.
We say that $f$ is {\em poly-time computable} if $T=T(\phi,n)$  computes $f(x)$ for each $x\in D$ and oracles $\phi$ of $x$, with running time bounded by a polynomial in $n$. 
\end{definition}

\subsection{Polynomial dynamics and invariant sets}

We let $p:\bC\to\bC$ be a complex polynomial with degree $d$ at least $2$.

We recall: 
(i) the \textit{Julia set} of $p$, $J$ or $J_p$, is the topological boundary of the filled Julia set $K=K_p$, which is the set of points with bounded orbits;  and
(ii) $p$ is \textit{hyperbolic} if it is uniformly expanding on its Julia set, in some Riemannian metric. That is, there is some neighborhood $\cN$ of $J$ and some ${L}>1$ and $c>0$ such that for all $z\in \cN$, $|(p^n)'(z)| \geq c {L}^n |z|,$ for all $n\in \bN$ for some Riemannian metric uniformly equivalent to Euclidean. As a consequence, some iterate of $p$ is strictly expanding ($c=1$) in the Euclidean metric. We refer to \cite{CarlGam} for more details about the dynamics of polynomials. 

For a hyperbolic, complex polynomial, the Julia set together with the attracting periodic points form the \textit{chain recurrent set}, $\cR$, where $\cR=\cap_{\alpha>0} \cR(\alpha)$ and $\cR(\alpha)$ is the set of all $\alpha$-recurrent orbits $z$ of $p$, that is, the set of points admitting an $\alpha$-pseudo orbit starting and ending at $z$.

\subsection{Distortion lemmas}

Finally,  we make significant use of the following, non-normalized version of the Koebe distortion theorem. Let $B^e(z,r)$  denote a Euclidean ball about $z$ of radius $r$, $|z|_e$ denote the Euclidean norm of $z$ in $\CC$, and $d_e(z,w)$ denote the Euclidean distance between $z$ and $w$ in $\CC$. 

\begin{theorem}[\cite{Geyer-metrics}]
\label{thm:koebe}
If $g\colon B^e(z_0,r) \to \CC$ is univalent, then 
\begin{equation} \label{eqn:koebe}
    \frac{d_e(z,z_0)}{\left( 1 + \frac{d_e(z,z_0)}{r} \right)^2} \leq 
    \frac{\left| {d_e(g(z),g(z_0))} \right |_e}{|g'(z_0)|_e} 
    \leq  \frac{d_e(z,z_0)}{\left( 1 - \frac{d_e(z,z_0)}{r} \right)^2}, 
\end{equation}
for any $z\in B^e(z_0,r)$.  
\end{theorem}

Note if $g$ is univalent on $B^e(z_0,r)$ then so are its iterates. In this article we leverage Equation~\ref{eqn:koebe} for the iterates $g^k$ of a univalent map $g$. 

We apply this theorem to develop the following bounds for the derivative of the map for points near $z_0$, which we employ in the proof of our main theorem. Since this is a general result we state it in this section. We only apply the lower bound in this article, but provide the upper bound as it is a general statement.  

\begin{proposition}
    
\label{prop:koebe-other-point}
Let $g\colon B^e(z_0,r) \to \CC$ be univalent and $0<a<r$. Set 
\begin{equation}
\label{eqn:gamma_r(a)}
\gamma=\gamma_r(a): = \max \left( 1-\frac{1}{(1+\frac{a}{r})^2}, \frac{1}{(1-\frac{a}{r})^2}-1 \right).    
\end{equation}

Then for any $w_0\in B^e(z_0,a),$ we have \begin{equation} \label{eqn:koebe-shift}
{|g'(z_0)|_e}
\left( \frac{1-3\gamma}{1+\gamma}
\right)
< {|g'(w_0)|_e} 
< {|g'(z_0)|_e}
\left(\frac{3(1+\gamma)}{1-\gamma} \right).
\end{equation}
\end{proposition}

\begin{proof}
From the definition of $\gamma$ and the distortion theorem (Theorem~\ref{thm:koebe}), for any $z\in B^e(z_0,a)$, we have 
\begin{equation}
\label{eqn:koebe-gamma}
d_e(z,z_0) (1-\gamma)
\leq   \frac{{d_e(g(z),g(z_0))}}{|g'(z_0)|_e} 
\leq  d_e(z,z_0) (1+\gamma).    
\end{equation}

Now, since $w_0\in B^e(z_0,a)$, let $\sigma>0$ be any constant small enough that $B^e(w_0,\sigma)\subset B^e(z_0,a)$.
For ease of notation we let $\zeta \in (0,a)$ satisfy $d_e(w_0,z_0)<\zeta$; so, $w_0\in B^e(z_0,\zeta)$ and $\sigma + \zeta < a$. 

Consider the point $w_1 \in \partial B^e(w_0, \sigma)$ which is the farthest point in $B^e(w_0,\sigma)$ from $z_0$. So, if you imagine the line between $z_0$ and $w_0$, the point $w_1$ is on the opposite side of $w_0$ from $z_0$. 
Thus we have $d(z_0,w_0)+d(w_0,w_1)=d(z_0,w_1)$, $d(w_0,w_1)=\sigma$, $d(z_0,w_0)<\zeta$, and $d(z_0,w_1) < \sigma + \zeta < a$.

Since both $w_0$ and $w_1$ lie in $B^e(z_0,a)$, 
 we can substitute either one in for $z$ or $z_0$ in Equation~\ref{eqn:koebe-gamma}.
We combine that with the triangle inequality to obtain the 
upper and lower bounds of 
the lemma.

First, for the upper bound we need only the upper bound of Equation~\ref{eqn:koebe-gamma} to get:
\begin{eqnarray*}
d_e(g(w_0),g(w_1)) 
& \leq & d_e(g(w_0),g(z_0))+d_e(g(z_0),g(w_1)) \\
 & \leq & |g'(z_0)|_e (1+\gamma)(d_e(w_0,z_0)
+d_e(z_0,w_1)) \\
 & = & |g'(z_0)|_e (1+\gamma)(2d_e(w_0,z_0)
+d_e(w_0,w_1)) \\
 & < &|g'(z_0)|_e (1+\gamma)(2\zeta+d_e(w_0,w_1)).    
\end{eqnarray*}

For the lower bound, we use both upper and lower bounds on Equation~\ref{eqn:koebe-gamma}. 
\begin{eqnarray*}
d_e(g(w_0),g(w_1)) & \geq & 
 d_e(g(z_0),g(w_1)) - d_e(g(z_0),g(w_0))     
\\ 
& \geq & 
|g'(z_0)|_e 
(d_e(w_1,z_0)(1-\gamma) - d_e(z_0,w_0)(1+\gamma)) 
\\
& = & |g'(z_0)|_e 
( d_e(w_1,z_0)-d_e(z_0,w_0) 
-\gamma d_e(w_1,z_0) -\gamma d_e(z_0,w_0) )
\\
& = & |g'(z_0)|_e  (d_e(w_0,w_1)-\gamma (d_e(w_0,w_1) + 2d_e(w_0,z_0))  )\\
& > & |g'(z_0)|_e  (d_e(w_0,w_1)(1-\gamma) - 2\gamma \zeta  )
\end{eqnarray*}
So, we have
\begin{equation*}
d_e(w_0,w_1)(1-\gamma) - 2\gamma\zeta
<   \frac{{d_e(g(w_0),g(w_1))}}{|g'(z_0)|_e}
<      (1+\gamma)(2\zeta+d_e(w_0,w_1)).
\end{equation*}

Now, we can apply the distortion theorem with $w_0$ as a center point to get:
$$d_e(w_1,w_0) (1-\gamma)
\leq \frac{{d_e(g(w_1),g(w_0))}}{|g'(w_0)|_e}
\leq   d_e(w_1,w_0) (1+\gamma).$$

Combining these yields both
$${|g'(z_0)|_e}(d_e(w_0,w_1)(1-\gamma) - 2\gamma\zeta ) \leq {{d_e(g(w_1),g(w_0))}}
\leq  {|g'(w_0)|_e} d_e(w_1,w_0) (1+\gamma) $$
and
$${|g'(w_0)|_e} d_e(w_1,w_0) (1-\gamma)
\leq
{{d_e(g(w_1),g(w_0))}}
\leq    {|g'(z_0)|_e}(  (1+\gamma)(2\zeta+d_e(w_0,w_1))).
$$
Ignoring the middle term and rearranging produces:
$$
{|g'(z_0)|_e}
\left( \frac{(1-\gamma)}{(1+\gamma) }
-\frac{ 2\gamma}{(1+\gamma)} \frac{\zeta}{\sigma}\right)
\leq  {|g'(w_0)|_e} $$
and
$${|g'(w_0)|_e} 
\leq  {|g'(z_0)|_e}\left(\frac{1+\gamma}{1-\gamma} \right) \left(2\frac{\zeta}{\sigma} +1 \right).$$
Combining these statements, and leveraging that this works for any $\sigma \in (0,\zeta),$ so the inequalities hold for $\frac{\zeta}{\sigma}=1,$ implies precisely Equation~\ref{eqn:koebe-shift}.
\end{proof}

Combining the proposition above with the distortion theorem allows us to bound, in terms of the derivative at a base-point, the size of the image of a ball near the basement, but with a different center. We also use this in our main theorem. 

\begin{lemma}
    \label{lem:proof-parts-1}
Let $g:B^e(z_0,r)\to\CC$ be univalent, $0<a<r$, and $\gamma=\gamma_r(a)$ given by Equation~\ref{eqn:gamma_r(a)}. 
Set $$c(\gamma) := \frac{(1-3\gamma)(1-\gamma)}{({1+\gamma})}.$$

Then for any $w_0\in B^e(z_0,a),$ and any $\sigma>0$ such that $B^e(w_0,\sigma)\subset B^e(z_0,a)$, we have: 
\begin{equation}
\label{eqn:lower-bound-W'}
\frac{1}{2}\text{diam}(g(B^e(w_0,\sigma))\geq |g'(z_0)|_e \cdot \sigma \cdot c(\gamma).
\end{equation}
\end{lemma}

\begin{proof}
Proposition~\ref{prop:koebe-other-point} provides:
$$
|g'(w_0)|_e \geq \left(\frac{1-3\gamma}{1+\gamma} \right) |g'(z_0)|_e.
$$ 
Consider any $w\in \partial B^e(w_0,\sigma)$.
Applying the distortion theorem around $w_0$ produces:
$$
d_e(g(w_0),g(w_0)) \geq |g'(w_0)|_e d_e(w_0,w)(1-\gamma).
$$
Combining the prior two inequalities and using the definition of $c(\gamma)$ and $\sigma$ yields
\begin{eqnarray*}
d_e(g(w_0),g(w)) 
&\geq &
d_e(w_0,w)\left(\frac{(1-\gamma)(1-3\gamma)}{1+\gamma} \right) |g'(z_0)|_e \\
& = & d_e(w_0,w)\ c(\gamma)\ |g'(z_0)|_e \\
&=& \sigma \ c(\gamma) \ |g'(z_0)|_e,
\end{eqnarray*}
so we have Equation~\ref{eqn:lower-bound-W'}.
\end{proof}

\section{Results}
\label{sec:results}

Let $p:\bC\to \bC$ be a complex polynomial with degree $d\geq 2$ and let $J_p$ denote the Julia set of $p$.
Given $N\in \bN$, we want to define a function $h_{J_p}(N,\cdot)$, on all ideal points with denominators $2^{N+2}$, which satisfies
Definition~\ref{defn:pic-according-to-h_C} 
for $C=J_p$. 

Moreover, the goal is that the value of $h_{J_p}(N,\cdot)$ can be determined (for all $N$ and all input numbers which are dyadic rationals with denominators $2^{N+2}$) via a Turing machine (a computer program for our purposes) whose running time  is bounded by a polynomial function in $N$. 

We briefly recall the fundamental difference between our algorithm and Braverman’s approach in \cite{Braverman2005} (see the discussion after Theorem~\ref{thm:J-computable} for more details). The algorithm \cite{Braverman2005} relies on complex-analytic tools—most notably the Poincar\'e metric and Pick’s theorem—to compute a neighborhood of the Julia set in which points outside the Julia set escape exponentially fast in the Poincaré metric. This escape rate is then algorithmically identified and subsequently translated into escape in the Euclidean metric. Such an approach inherently depends on the holomorphic structure and does not extend to non-holomorphic systems.

In contrast, our method is based entirely on dynamical tools, specifically quantitative shadowing and chain recurrence. We compute a collection of small reference boxes that cover the Julia set, each of which has nonempty intersection with the Julia set. Moreover, we determine a constant 
$L>1$ such that a suitable (computable) iterate of the polynomial has derivative greater than 
$L$ in the Euclidean norm throughout each box. Using these expansion bounds together with quantitative escape and distortion estimates, we obtain a criterion to decide whether an arbitrarily small ideal ball contains a point of the Julia set.

We note that our algorithm for poly-time computability of $J_p$ for a hyperbolic $p$ has a 
running time which is a constant plus a linear function of $N$. This constant will be an initial ``overhead'' computation coming from applying parts of the algorithms of \cite{BoydWolf-Skew1}, adapted to a polynomial $p$. 
The algorithms and results of \cite{BoydWolf-Skew1} for polynomial skew products $f(z,w) = (p(z), q(z,w))$ of $\CC^2$ can be straightforwardly modified to apply in the simpler setting of polynomial maps $p(z)$ of $\CC$. In the one-dimensional setting, there are no saddle sets, only the expanding Julia set $J_p$ and a finite number of attracting periodic orbits, $A_p$. 

\subsection{Overhead/first step computation}
\label{sec:firststep}

In this subsection, we build the ``first step'' of our poly-time algorithm. This step is independent of the desired precision ($2^N$-approximation of $J_p$). In other words, once the polynomial is identified, this part of the algorithm runs once, and halts if $p$ is hyperbolic. Its results can then be used for the $2^N$-approximation of the corresponding Julia set, for any $N$. For readability, we split this initial algorithm in several parts. First, we show how to build a collection of boxes containing $J$ on which the map is expanding.

\begin{lemma}
\label{lem:make_N_1}
    There exists a Turing machine which based on input of an oracle of a  hyperbolic polynomial $p:\bC\to\bC$  of degree $d\geq 2$, computes positive dyadic rationals $\delta$ and $\lambda'$, a positive integer $\nu$, and a finite collection of closed square boxes $\sB'_n= \{B_j\}$ (with disjoint interiors, which are a subset of the boxes in a $2^n \times 2^n$ grid on a large box $[-R,R]^2$ containing $\cR(p)$), such that:
    \\ (i) the boxes' union contains $J$, i.e., $\cB'_n = \cup_j B_j \supset J$, and (ii) $p^\nu$ is expanding on the boxes' union; in particular,
 $|Dp^{\nu}(z)|>1+\lambda'$ for all $z\in \cB'_n$.
 
 (iii) Further, on $\cN_1$ defined as the collection of boxes with the same center points as $\cB'_n$ but with radii increased by $2\delta$ (i.e., $\cN_1:=\bigcup_{ B\in \sB'_n } \cN(B,2\delta) = \cN(\cB'_n,2\delta)$, thus $\cN_1 \supset \cB_n' \supset J$), we have $p^\nu$ is still expanding, specifically:
$ |Dp^{\nu}(z)|>1+\lambda'/2$ for all $z \in \cN_1$.
\end{lemma}
\begin{proof}
We prove by providing the algorithm. 
First, Algorithm 3.2 of \cite{BoydWolf-Skew1} can be applied to $p$ to produce a ``box chain model'' which includes a set $\cB_n$ which is the union of a collection of boxes, $\sB_n = \{B_j\} $, with the boxes all of the same side-length $\ep_n$, with pairwise disjoint interiors;
in fact, these boxes are a subset of the collection of boxes formed from a $2^n \times 2^n$ grid on a large box $[-R,R]^2$, where $R$ is a power of $2$ and is calculated to be guaranteed to contain the chain recurrent set $\cR$ of $p$, and hence also contains $J=J_p$ and the attracting periodic points $A_p$ of $p$. The ``model'' is this collection of boxes together with a graph $\Gamma_n$ whose vertices are the boxes $B_j$ and there is guaranteed to be an edge from $B_k$ to $B_j$ if $p(B_k)$ intersects $B_j$ (and there is guaranteed not to be an edge if $p(B_k)$ is some given distance away from $B_j$); moreover, every box in $\Gamma_n$ lies in a cycle. Further, decomposing $\Gamma_n$ into strongly connected components yields ``box chain components'' $\cB_{n,i}$ such that if the boxes are sufficiently small, $J$ and each attracting periodic orbit are contained in different box chain components: there are no edges in $\Gamma$ connecting different box chain components.

Algorithm 3.8  (which builds on Algorithm 3.2) of \cite{BoydWolf-Skew1}, can be adapted to the polynomial map $p$ to compute a sufficiently large $n\in \bN$  and a positive integer $\nu$ such that the box chain model $\cB_n$ consists entirely of box chain components $\cB_{n,i}$ with  each determined to be of ``expanding'' or ``contracting'' type (for polynomials there is no saddle type), and this algorithm halts if $p$ is hyperbolic; specifically, there exist $\lambda_i>0$ such that either (a) $|Dp^\nu(z)| > 1+\lambda_i$ for all $z$ in all boxes in $\cB_{n,i}$, or (b) $|Dp^\nu(z)| < 1-\lambda_i$ for all $z$ in all boxes in $\cB_{n,i}$. 

We let $\sB'_n$ denote the collection of all boxes determined to be of expanding type, and we let  $\cB'_n$ be their union. We have $|Dp^\nu(z)| > 1+\lambda',$ $\forall z \in \cB_{n}'$ (where $\lambda'$ is the minimum of the $\lambda_i$'s for all expanding-type components $\cB_{n,i}$). Hence we know if $p$ is hyperbolic, $\cB'_n \supset J$.

By applying an argument involving the computable modulus of continuity of the derivative of $p^\nu$, as was done similarly in \cite{BoydWolf-Skew1}, we can compute a bound for a larger neighborhood of $\cB'_n$ on which expansion is only $1+\lambda'/2$. In particular,  we can compute a lower bound $\delta>0$, and make it a bit smaller so that $\delta$ is a dyadic rational, such that if we replace each box $B$ in $\sB'_n$ by the slightly larger box which is its $2\delta$-neighborhood $\cN(B,2\delta)$ in our $L^\infty$ metric (using $2\delta$ here is convenient as later, $\delta$ will be significant), we have a collection of boxes with the same center points but with side-lengths larger by $4\delta$ (so the larger box side-length is $4\delta+\ep_n$ where $\ep_n$ is the box side-length of $\cB'_n$), such that
 all points in these larger boxes satisfy an inequality of expansion but only with a guarantee of weaker expansion $L:=1 + \lambda'/{2}$. That is, $|Dp^\nu(z)| > 1+\lambda'/2=L$ for all $z \in \cN(\cB'_n,2\delta)$.  These larger boxes may overlap, but that is no problem.

 Set $\cN_1:= \bigcup_{B_j\in \sB'_n} \cN(B_j,2\delta)$, so the Hausdorff distance of $\cB_n'$ and $\cN_1$ is exactly $2\delta$ (i.e., $d(H(\cB'_n,\cN_1)=2\delta$). Since $J\subset \cB'_n$, it follows that $J\subset \cN_1$, which completes the proof of the lemma. Note we also have $d_H(J,\cN_1)\geq 2\delta$ and $\cN_1\supset \cN(J,2\delta).$ 
 \end{proof}

Next we build onto Lemma~\ref{lem:make_N_1}. We use shadowing to refine the boxes of $\cN_1$ into a smaller set still containing $J$ but for which each box contains a point of $J$. 

\begin{proposition} \label{prop:firststep}
    Given a hyperbolic polynomial $p:\CC\to \CC$ of degree $d\geq 2$, suppose $L>1, \delta>0$, $\nu\in \bN^+$  and  $\cN_1$ are produced by the Turing machine in Lemma~\ref{lem:make_N_1}. So, we have
  $|Dp^{\nu}(z)|>L$ for all $ z \in \cN_1 \supset \cN(J,2\delta).$
  
  Then there exists a Turing machine which on input of an oracle of a real number $b>0$ computes a set which is a finite union of closed square boxes, $\cN_2=\cup B_j$,
     with each box of side-length $\beta'\in (0,b]$, such that
\begin{itemize}
    \item [(i)]
    the union of these boxes contains the Julia set:  $J_p \subset \cN_2$, 
    \item [(ii)]  each box $B_j$ which is an element of the union of boxes defining $\cN_2$ contains a point of $J_p$:  $B_j \cap J_p \neq \emptyset$ 
    for all $B_j$ in  $\cup B_j = \cN_2$,
      \item[(iii)]  $p^\nu$ has the same expansion on $\cN_2$ as on $\cN_1$: $|Dp^{\nu}(z)|>L,$ $\forall z \in \cN_2$, and 
      \item[(iv)]  the Hausdorff distance from $\cN_2$ to $\cN_1$ is at least
      $\delta$, and a $\delta$-neighborhood of $\cN_2$ is contained in $\cN_1$: $d_H(\cN_2,\cN_1) > \delta$ and $\cN(\cN_2, \delta) \subset \cN_1$. 
\end{itemize}
We call these boxes ``reference boxes''. 
\end{proposition}

\begin{proof}
    To prove this lemma, we supply the algorithm. The first step is to invoke Lemma~\ref{lem:make_N_1}, creating the sets $\cB'_n, \cN_1$, and finding the constants $L,\lambda', \delta$ and $\nu$. Recall from the proof of the lemma, $\cB'_n$ is a level $n$ box chain model of $p$,  and we have $|Dp^{\nu}(z)|>L=1+\lambda'/2$ for all $ z \in \cN_1 \supset \cN(J,2\delta).$

The boxes are nested upon refinement: $\cB_{n} \subset \cB_{n+1}$, you may subdivide the boxes of $\cB_n$ and repeat the algorithm to obtain sub-boxes which are in $\cB_{n+1}$.

Further, Theorem 3.6 of \cite{BoydWolf-Skew1} includes (as a consequence of Algorithm 3.2 of \cite{BoydWolf-Skew1}) that there is a computable sequence $\alpha_n\downarrow 0$ as $n\to \infty$ such that $\cR \subset \cB_n \subset \cR(\alpha_n)$. Hence, the boxes defining $\cB_n$ are contained in $\alpha_n$-pseudo-orbits, and an upper bound for $\alpha_n$ is computable and converges to $0$ as $n\to\infty$.  
Also, the box side-length $\ep_n$ (for boxes defining $\cB'_n$) satisfies $\ep_n < \alpha_n$.

We now use a modified part of Algorithm 3.10 of \cite{BoydWolf-Skew1}--changed slightly more than just adapting to the one-dimensional case. Here, we describe the outline for this step. We focus on the ``(a) Type $J_2$-algorithm'' part of Algorithm 3.10, which locates a neighborhood of the Julia set (closure of repelling periodic points) to a given precision. Rather than applying this to yield the precision desired based on $N$,  what we need is only to subdivide the boxes in $\cB'_n$ just enough to make smaller boxes such that given the expansion $1+\lambda'$, the pseudo-orbits in these smaller boxes are  $\beta$-shadowed by true orbits for some $\beta>0$. Then for boxes with radius $\beta$ larger, these slightly larger boxes with the same center point each contains a point of $J$, a true orbit; so the sidelength of these boxes which each contain a point of $J$ is our $\beta'$.
Finally, we need the $\beta'$-boxes to still be in the neighborhood of $J$ in which $p$ is expanding (that is, $\cN_1$). This creates another upper bound on $\beta'$ (as does the provided $b$). The union of these $\beta'$-sidelength boxes defines our set $\cN_2$.  
 We provide more details in the appendix, establishing items (i)--(iv). 
\end{proof}
Combining these results produces nested sets  with
$J \subset 
\cN_2 \subset \cN_1 $
and the following properties for $g=p^\nu$:
\begin{itemize}
    \item $|Dg|>1+\lambda'$ on $J$, 
    \item $|Dg|>1+\lambda'/2=L$ on both $\cN_2$ and $\cN_1$,
\item 
    $d_H(\cN_1,J)\geq 2\delta$ (since $d_H(\cB'_n,\cN_1)=2\delta$ and $J\subset \cB'_n$)  and $\cN_1\supset \cN(J,2\delta),$ 
\item $d_H(J,\cN_2) < \beta'$, where $\beta'$ is the side-lengths of the boxes in $\cN_2$, and
\item $d_H(\cN_2,\cN_1)> \delta$
 and $\cN(\cN_2, \delta) \subset \cN_1$.
\end{itemize}

\subsection{Poly-time computability of $J$}
\label{sec:mainresult}

In this subsection, we first provide our algorithm which incorporates the steps of the prior subsection to create the approximation of $J$ to the desired precision. Then we prove the algorithm as claimed provides the desired approximation in linear time for hyperbolic polynomials.

\begin{algorithm} \label{alg:computeJ}
\textbf{Computing a $2^N$-approximation of $J$.} Consider a hyperbolic polynomial $p:\CC\to\CC$  of degree $d\geq 2$ which is given by an oracle of its coefficients. Our first step is as follows.

\textbf{Step 1(a).}
First, if $p$ is hyperbolic, we can build a collection of boxes on which an iterate $p^\nu$ is expanding. That is, we invoke Lemma~\ref{lem:make_N_1},  which succeeds if $p$ is hyperbolic, to
compute a positive integer $\nu$, dyadic rationals $\delta>0, L>1$,
 and the set $\cN_1$ which is the union of a collection of boxes on which $|Dp^\nu(z)|>L$.  
If $p$ is not hyperbolic, this step will never halt. If this step halts, then the algorithm has just proven that $p$ is hyperbolic (see Lemma 3.9 of~\cite{BoydWolf-Skew1}). 
 
 For ease of notation, set $g=p^\nu$. Note since $|g'(z)|$ is bounded below by $L>1$, the map $g$ is conformal on $\cN_1$.

\textbf{Step 1(b).}
Next, we want to invoke Proposition~\ref{prop:firststep}  to construct our ``reference boxes'' defining $\cN_2$. To do so, recall the proposition allows us to
 input a $b>0$ that is an upper bound for the side-length $\beta'$ of the boxes defining $\cN_2$, though, we never said why such a $b$ might be needed. Now, we provide the calculation defining the $b$ we need.

First set $r=r(\delta):=3\delta/4$ (which is a bit arbitrary but this choice is convenient for calculations) and choose any $\gamma \leq 0.1$ (we see  why that $\gamma$ is useful in the proof of (iii) in Lemma~\ref{lem:step2-conditions}). Then use Equation~\ref{eqn:gamma_r(a)} defining $\gamma_r(a)$ in terms of $a$ and $r$ to determine $a$ from this $r$ and $\gamma$. Now, set $c=c(\gamma) = \dfrac{(1-3\gamma)(1-\gamma)}{({1+\gamma})}$ (matching Lemma~\ref{lem:proof-parts-1}), and using $\gamma$, $c(\gamma)$, and $a$, set:
\begin{equation}
    \label{eqn:defn-b}
b:=\min\left\{ \frac{\delta \ c(\gamma)}{16\sqrt{2}(1+\gamma)}, a\right\}.
\end{equation}
The reason this value of $b$ is needed is demonstrated in the proof of Theorem~\ref{thm:J-computable} (more precisely, this $b$ is used at the end of the proof of Lemma~\ref{lem:step2-conditions}).

With $b$ is determined, we now invoke Proposition~\ref{prop:firststep} to compute the union  $\cN_2$ of boxes of computable side-length at most $\beta' \in (0 , b)$, such that 
 (i) $J \subset \cN_2$, (ii) each box in the union defining $\cN_2$ contains a point of $J$, (iii) $|Dg(z)|>L>1$ for every  $z\in \cN_2$, and  
(iv) $d_H(\cN_2,\cN_1)> \delta$
 and $\cN(\cN_2, \delta) \subset \cN_1$, where the boxes defining $\cN_1$ and $\delta>0$ were constructed in Step 1(a) as a set on which the iterate $g=p^\nu$ is expanding.

\textbf{Step 1(c).} 
We explain how to draw the $2^N$-picture of $J$. (Why it actually is a $2^N$ picture is not explained here, but proven  following the algorithm.)

First, recall (at the beginning of the algorithm described in Lemma~\ref{lem:make_N_1}) we calculated a large box $[-R,R]^2$ containing $J$.
So, we can define $h_{J}(N,z')=0$ for all $z'\in \CC \setminus [-R-1,R+1]^2$.
Now we consider only $z' \in [-R-1,R+1]^2$.

In the following, it turns out we need $N\geq N'$, where $N'$ is the smallest positive integer which satisfies both of the following::
\begin{equation}
2^{-N'-1}<a \ \text{ and }  \ 2^{-N'}<\delta, \label{eqn:define_N'}
\end{equation}
for the $\delta>0$ determined in Step 1(a) and the $a$ chosen in Step 1(b).
 
To draw a $2^N$-picture of $J$ for some $N<N'$ we simply draw the more precise  $2^{N'}$-approximation. 

Now consider any $N\geq N'$ (so Equation~\ref{eqn:define_N'} holds for $N$). We perform the following for every ideal point $z'\in [-R-1,R+1]^2$, a dyadic rational with coordinate denominators in $2^{N+2}$.

 Recall that $B(z', r)$ denotes the ball (box in our norm) about $z'$ of radius $r$. 
 For ease of notation, set $\cU':=\overline{B(z', 2^{-N-2})}$ and $2\cU':= \overline{B(z', 2^{-N-1})}$.
 
\textbf{Step 1(d).} This is a ``base case''. Test whether $\cN_2\cap  \cU' = \emptyset$. 
If so, stop and report $0$ for $z'$ (i.e., set $h_J(N,z')=0$ and halt the algorithm for $z'$; go to the next ideal point to consider). 

\textbf{Step 1(e).} This is a  ``base case''. Test whether there is a reference box of $\cN_2$ entirely contained in $2\cU'$. 
If so, stop and report $1$ for $z'$ (set $h_J(N,z')=1$). 

\textbf{Step 1(f).} This is the final ``base case''. Test whether $2^{-N-2} > \beta'/2.$ 
If so, stop and report $0$ for $z'$.

    \medskip
    
\textbf{Step 2.}
If this step is reached, then $\cN_2\cap \cU'\neq \emptyset$,
but there is no reference box entirely contained in $2\cU'$.
Also, $2^{-N-2} \leq  \beta'\sqrt{2}.$
We need to determine whether to report $0$ or $1$ for $z'$.

First we provide an outline of this step. We leverage the distortion theorem (Theorem~\ref{thm:koebe}),  for $z_0=z'$, and iterates of $g=p^\nu$. To do so we must replace $\cU'$ and $2\cU'$ with appropriate Euclidean balls centered at $z'$. 

 The idea is that for iterates $g^k$, the distortion theorem provides over and under-estimates for the images under $g^k$ of these ball approximations to $\cU'$ and $2\cU'$, in terms of the derivatives $|Dg^k(z')|_e$. Since for large $N$, the sets $\cU'$ and $2\cU'$ are small, the radii of these balls are very small compared to the size of a ball of conformality, 
 so the over- and underestimate images do not differ by much. 
 
 If the ideal point $z'$ is close to $J$, there is a point of the Julia set near $z'$ which guarantees expansion along the orbit of $z'$, for a while at least. Consequently, the diameter of the image of $\cU'$ grows as we iterate, and 
 either the balls we are tracking about the orbit of $z'$ move  far enough away from our neighborhood of expansion in which case we know that $z'$ wasn't sufficiently close to $J$ and we report $0$, or we can iterate sufficiently many times to guarantee that the diameter of the images under $g^k$ of $\cU'$  is larger than the radius of the reference boxes in $\cN_2$, thus guaranteeing a point of $J$ is in one of these images, hence in the original box by invariance of $J$. In other words, we analyze if there is either a reference box entirely contained inside the underestimate for the ball approximation to $g^k(2\cU')$, which happens if there was a point of $J$ in $g^k(\cU')$, hence in $\cU'$, in which case we report $1$, \textit{or} we analyze if all reference boxes are outside of the over-estimate ball approximation to $g^k(\cU')$, which happens if $J$ is outside of $g^k(\cU')$, so by invariance $\cU'\cap J=\emptyset$ and we report a $0$. 
We must also apply the distortion to bound the underestimate of the larger ball sufficiently away from the overestimate for the smaller ball.  
 Finally, we calculate the number of iterates we need to test this procedure, and report $0$ if neither test criteria is reached after this number of iterates.

Now we present the detailed algorithm for this step. 

\textbf{Step 2(a).} 
We first calculate the maximum number of iterates that have to be considered. (We use this in the proof of Lemma~\ref{lem:Step2f-correct}.)
Define
\begin{equation}
\label{eqn:defn-k-N}
k_N:= \left \lceil \frac{1}{\log_2(L)}
\left( 
N + 3 + \log \left(\frac{2\beta'\sqrt{2}}{c(\gamma)}\right)
\right) \right\rceil,
\end{equation}
where $\lceil \cdot \rceil$ is the ceiling function, and recall from Step 1(b), $\gamma\in (0,0.1]$ and  $c=c(\gamma) = (1-3\gamma)(1-\gamma)/({1+\gamma})$, $L$ is the lower bound on expansion of $|Dg|$ in $\cN_1$ (and $\cN_2$), and $\beta'$ is the side-length of the reference boxes defining $\cN_2$. 
It is important that $k_N$ is linear in $N$.

We need balls rather than boxes so as to apply the distortion theorem. For $\sigma_N := 2^{-N-3},$ we consider  
\begin{equation} \label{eqn:Omega-defn}
\Omega_U :=   B^e(z', 4  \sigma_N)     \text{ and } \Omega_L:= B^e(z',3  \sigma_N). 
\end{equation}
Thus, the width of the annulus $\Omega_U\setminus \Omega_L$ is $\sigma_N$. 
An elementary calculation confirms that $\cU' \subset \Omega_L \subset \Omega_U\subset 2\cU'$, see Figure~\ref{fig:ideal-ball-sandwich}.
\begin{figure} \centering
\includegraphics[width=.457\textwidth]{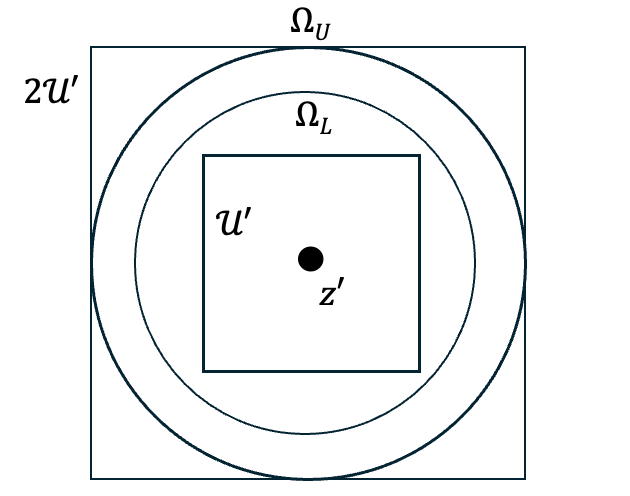}
\includegraphics[width=.45\textwidth]{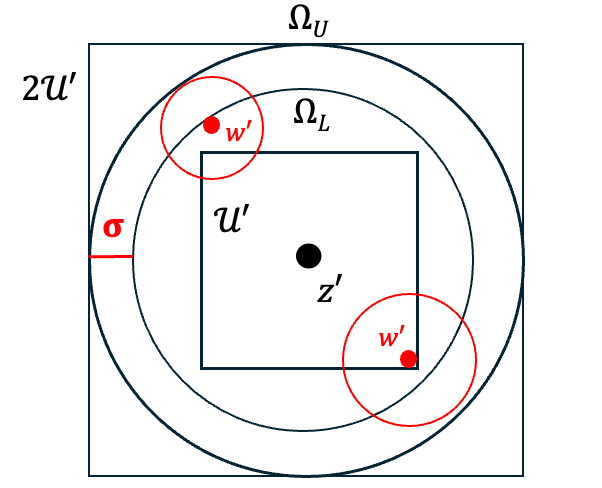}
\caption{\label{fig:ideal-ball-sandwich}
Left: In Algorithm~\ref{alg:computeJ}, we consider boxes  $\cU'$ and their ``doubled'' boxes with the same center an ideal point $z'$ but double diameter $2\cU'$. We also consider Euclidean balls nested between these, so $\cU' \subset \Omega_L \subset  \Omega_U \subset 2\cU'$. Right: Additionally, we consider (hypothetical) smaller balls inside of $\Omega_U$, with center points in $\Omega_L$.}
\end{figure}

\medskip

\noindent Next, consider $k=1$ up to (at most) $k=k_N$, in turn.  
For each such $k$: 

\textbf{Step 2(b).} 
Calculate $g^k(z')$ and $|Dg^k(z')|_e$.

\textbf{Step 2(c).} 
Test  whether 
the ball 
$B^e(g^k(z'), |Dg^k(z')|_e (3\sigma_N) (1+\gamma))$ (which we justify in the proof of the algorithm's correctness is the overestimate for $g^k(\Omega_L)$ if we get to this step for $k$) is disjoint from $\cN_2$.
If so,  stop and report $0$ for $z'$. 

\textbf{Step 2(d).} Test if  
$B^e(g^k(z'), |Dg^k(z')|_e (4\sigma_N) (1-\gamma))$
(which we show in the proof below is the underestimate of $g^k(\Omega_U)$)
entirely contains any reference box of $\cN_2$. 
If this holds then we stop the algorithm for $z'$ and report $1$. 

\textbf{Step 2(e).} 
Test if $|Dg^k(z')| \sigma_N c(\gamma) > 2\beta'\sqrt{2}.$ If so, stop at this $k$ and report $0$ for $z'$. 
(We explain in the proof below what this quantity represents).

If we pass through Steps 2(b,c,d,e) and do not report $0$ or $1$ for a $k<k_N$, we proceed to repeat the loop of Steps 2(b,c,d,e) with $k$ replaced by  $k+1$. If we get  past the final loop for $k=k_N$, we move on to Step 2(f). 

\medskip

\textbf{Step 2(f).} 
Finally, if we get through all $k\in \{1, \ldots, k_N\}$ and have not yet stopped the algorithm and reported $0$ or $1$ for $z'$, we report $0$ for $z'$. 
\qed
\end{algorithm}

\medskip

We prove Theorem~\ref{thm:J-computable} by first going through the steps of the algorithm and showing that it produces a $2^{N}$-approximation of $J$. Then, we prove that the algorithm halts. Finally, we examine the running time. 

First we examine correctness of Step 1.

\begin{lemma}\label{lem:step1-correct}
    Let $p$ be a polynomial map of $\CC$ of degree $d\geq 2$. Beginning to apply Algorithm~\ref{alg:computeJ}, suppose we are examining $z'$ an ideal point of denominator $2^{N+2}$ for an $N\geq N'$.  If Algorithm~\ref{alg:computeJ} halts during one of the Step 1 tests (at Steps 1(d, e, or f)) and reports $0$ or $1$ for $z'$, it does so consistently with the definition of $2^N$-approximation of $J$. That is,
    \begin{itemize}
     \item if \textbf{Step 1} reports $0$ for $z'$, then $d(z',J)>2^{-N-2}$;
     
     \item if \textbf{Step 1} reports $1$ for $z'$, then $d(z',J) \leq 2^{-N-1}$.
    \end{itemize}
\end{lemma}

\begin{proof}
We explained in Step 1(a) of Algorithm~\ref{alg:computeJ} (and in Section~\ref{sec:firststep}) that if Step 1(a) succeeds, then we have established that $p$ is hyperbolic, and built a union of boxes $\cN_1$ containing $J$ and on which $p$ is expanding, and a union of smaller boxes $\cN_2$ containing $J$ such that each box contains a point of $J$, and 
$\cN_1\supset \cN(\cN_2,\delta)$. (Note if $p$ is not hyperbolic, Step 1(a) never halts so we do not make it to the tests Step 1(d, e, or f) to report anything for any $z'$, so this is why we do not need to assume $p$ is hyperbolic for the lemma, that is implied by the assumption the algorithm makes a report in Steps 1(d, e, or f).) 

\textbf{Step 1(d)} contained the first criterion which we checked for an ideal point to determine $h_J(N,z')$ in Algorithm~\ref{alg:computeJ}; if $\cN_2\cap \cU'=\emptyset$ then report $0$. 
Since $J\subset \cN_2$, in that case we have $z'$ is at least $2^{-N-2}$ (the radius of $\cU'$) from $J$. We are allowed to report $0$ for any point which is $\geq 2^{-N-2}$ from $J$. 

\textbf{Step 1(e)} contained the second check. Since every reference box defining $\cN_2$ contains a point of $J$, if there is a reference box entirely contained in $2\cU'$,
then there is a point of $J$ in $2\cU' = \overline{B(z', 2^{-N-1})},$ so $d(J,z')\leq 2^{-N-1}$. So we are allowed to report $1$.

\textbf{Step 1(f)}. If we get here, we did not report $1$ in Step 1(e) which means $2\cU'$
did not entirely contain a reference box. But,  if there was a point of $J$ in $\cU'$, then since the (box) radius of the reference box is smaller than the difference between $2\cU'$ and $\cU'$,
the reference box of that point would have been entirely contained in $2\cU'$,
and we would have halted in the prior step. Since that did not happen, there is no point of $J$ in $\cU'$
so we may report $0$. 
\end{proof}

Next, for the correctness of Step 2, we must first establish some conditions that are true if we make it to Step 2 of the algorithm.

\begin{lemma} \label{lem:step2-conditions}
     Let $p$ be a polynomial map of $\CC$ of degree $d\geq 2$. Beginning to apply Algorithm~\ref{alg:computeJ}, suppose we are examining $z'$ an ideal point of denominator $2^{N+2}$ for an $N\geq N'$.   If Algorithm~\ref{alg:computeJ} makes it to Step 2, then we know:

    \begin{itemize}
\item [(i)] $B^e(z',r)\subset \cN_1$, in particular $g$ is conformal on $B^e(z',r)$, for $r=3\delta/4$.
  \item[(ii)] If $w'\in \Omega_L$, then $d(w',\Omega_U)\geq \sigma_N$, and for any $i\in \bN$,
{$$\frac{1}{2} {\rm diam}(g^{i}(B^e(w',\sigma_N))\geq |Dg^{i}(z')|_e \sigma_N c(\gamma).$$} 
\item[(iii)] Suppose for some $i\in \bN$, we know  $|Dg^i(z')|_e\sigma_N c(\gamma) > 2\beta'\sqrt{2}$.
If $\Omega_L$ contains a point of $J$, then there is a reference box of $\cN_2$ inside the underestimate of $g^i(\Omega_U)$ (the ball $B^e(g^i(z'), |Dg^i(z')|_e (4\sigma_N)(1-\gamma))$). 
 \item[(iv)]  For any $k\in  \{1,\ldots,k_N\},$  if we make it to the start of the $k$-loop of Step 2(b), then the all iterates $z', \ldots, g^{k-1}(z')$ are contained in $\cN(\cN_2,\delta/2)$.
\end{itemize}    
\end{lemma}

\begin{proof} 
    Note first that if we make it to Step 2, then Step 1(a) successfully computed the neighborhoods $\cN_2, \cN_1$ (rather than running forever).  As a consequence, it proves that the map $p$ is hyperbolic; thus this need not go into the assumptions. 
    
    But also, if we make it to Step 2, none of the tests in Step 1 triggered a report. So: $\cN_2\cap \cU'\neq \emptyset$ (though
 there is no reference box entirely contained in $2\cU'$, and since $2^{-N-2} \leq  \beta'\sqrt{2},$  the radius of $2\cU'$ is at most twice the Euclidean diameter of a reference box).

\underline{(i)}: First we show $B^e(z',r)\subset \cN_1$, in particular $g$ is conformal on $B^e(z',r)$, for $r=3\delta/4$. 
Since Step 1(d) did not triggered a report, we know $\cN_2\cap \cU'\neq \emptyset$. Thus $d(z',\cN_2)\leq 2^{-N-2}$.
Recall that $\cN_1\supset \cN(\cN_2,\delta).$ 
By Equation~\ref{eqn:define_N'}, we have $\delta > 2^{-N'}\geq 2^{-N}$, thus $2^{-N-2} <\delta/4,$ so $d(z',\cN_2)\leq \delta/4,$ hence $B^e(z',3\delta/4=r)\subset \cN_1$. Finally, recall that $g$ is conformal on $\cN_1$ since  we have $|Dg(z)| > L > 1 >0$ there.

\underline{(ii)}:
Next we show if $w'\in \Omega_L$, then $d(w',\Omega_U)\geq \sigma_N$, and for any $i\in \bN$,
{$\frac{1}{2} {\rm diam}(g^{i}(B^e(w',\sigma_N))\geq |Dg^{i}(z')|_e \sigma_N c(\gamma).$} 
We aim at applying Proposition~\ref{prop:koebe-other-point} and Lemma~\ref{lem:proof-parts-1} for $g^i=(p^\nu)^i,$ and $r, a, \gamma$ as defined in Step 1(b) ($ r=(3\delta/4)$, $\gamma \leq 0.1$ and  $a< r$ based on $\gamma_r(a)$ defined by Equation~\ref{eqn:gamma_r(a)}), and for $z_0=z'$, $w_0=w'$ any point of 
$\Omega_L$, and $\sigma=\sigma_N$($=2^{-N-3}$ as defined in Step 2(a)), see Figure~\ref{fig:ideal-ball-sandwich}, Right. 
Since we are in Step 2, we know $\cU'\cap \cN_2\neq \emptyset$ so  $z'\in \cN(\cN_2,\delta/2)\subset \cN_1$ and $\Omega_L\subset \Omega_U \subset  \cN_1$ as well since 
$2^{-N-1}<\delta/2$. So since by Equation~\ref{eqn:define_N'}, we chose $N'$ in Step 1(c) so that $2^{-N'-1}<a,$ we know $\Omega_L \subset B^e(z',a)$. Hence $w'\in B^e(z',a)$, and Lemma~\ref{lem:proof-parts-1} yields (ii).

\medskip

\underline{(iii)}:  Now we show that if for some $i\in \bN$ we have  $|Dg^i(z')|_e\sigma_N c(\gamma) > 2\beta'\sqrt{2}$, \textit{and} $\Omega_L$ contains a point of $J$, then there is a reference box of $\cN_2$ inside the underestimate of $g^i(\Omega_U)$ (the ball $B^e(g^i(z'), |Dg^i(z')|_e (4\sigma_N)(1-\gamma))$). 

Note since there is a point of $J$ in $\Omega_L$, say $w'$, then we know $z' \in \cN_1$, so  $B^e(g^i(z'), |Dg^i(z')|_e (4\sigma_N)(1-\gamma))$ is the underestimate for $g^i(\Omega_U)$. 

Consider the ball $W'=B^e(w',\sigma_N)$.
We showed above in (ii) that 
for $w'\in \Omega_L$ (hence $d(w',\Omega_U)\geq \sigma_N$), we know 
$\frac{1}{2}{\rm diam} (g^i(W')) \geq |Dg^i(z')|_e \sigma_N c(\gamma).$

Thus, this implies by our hypothesis that 
$\frac{1}{2}{\rm diam} (g^i(W'))> 2\beta'\sqrt{2}$, i.e.,
the radius of $g^i(W')$ is greater than twice the Euclidean diameter of the reference boxes.
So, for $w' \in J \cap \Omega_L$, we know $g^i(w') \in J\cap g^i(\Omega_L)$, hence given its size, $g^i(B(W'))$ must contain the reference box that contains the Julia set point $g^i(w')$;
in fact this reference box is at least $\beta'\sqrt{2}$ away from the boundary of $g^i(W')$. Also, by choice of $\sigma_N$, since $W' \subset \Omega_L$,  the reference box inside of $g^i(W')$ is inside of $g^i(\Omega_U)$, also a distance of at least $\beta'\sqrt{2}$ from the boundary $\partial g^i(\Omega_U)$. 

Observe that if the difference of the radii between the underestimate of $g^i(\Omega_U)$ and the overestimate of $g^i(\Omega_L)$ is at least $\beta'\sqrt{2},$ then we know that the reference box about $g^i(w')\in g^i(\Omega_L)$ is contained entirely inside of the underestimate of $g^i(\Omega_U)$. We show this is true. In fact, this difference in radii is at least:
{$$|Dg^i(z')|_e4\sigma_N(1-\gamma) - |Dg^i(z')|_e 3\sigma_N(1+\gamma) = |Dg^i(z')|_e\sigma_N(1-7\gamma).$$}
We set $c(\gamma)=(1-3\gamma)(1-\gamma)/(1+\gamma),$ so by an easy calculation,
by choice of $\gamma \leq 0.1$, we have  $(1-7\gamma) \geq c(\gamma)/2$.
(Simplifying this inequality leads to a quadratic which opens down, with two roots, one negative, and the positive root has solution $\gamma = -(8-\sqrt{132}/34) \approx 0.1026213$, so any $\gamma \in (0,0.1]$ satisfies this desired inequality). 

Thus, $|Dg^i(z')|_e\sigma_N(1-7\gamma) \geq |Dg^i(z')|_e\sigma_N c(\gamma)/2.$ Therefore, since by hypothesis this latter quantity is $> \beta'\sqrt{2}$, we know the 
 the minimum difference between the underestimate of $g^i(\Omega_U)$ to the boundary of the overestimate of $g^i(\Omega_L)$, is greater than the Euclidean diameter of a reference box, and we
conclude that the reference box about $g^i(w')$ is contained entirely inside of the underestimate of $g^i(\Omega_U)$.

 \medskip

\underline{(iv)}: Finally, we show for any $k\in  \{1,\ldots,k_N\},$  if we make it to the start of the $k$-loop of Step 2(b), then the iterates $z', \ldots, g^{k-1}(z')$ are contained in $\cN(\cN_2,\delta/2)$. We proceed by induction on $k$. 

For the base case $k=1$, if we make it to Step 2 for the first time, we know that $\cU' \cap \cN_2\neq \emptyset.$ Thus 
$g^0(z')=z'\in \cN(\cN_2,2^{-N-2}) \subset \cN(\cN_2,\delta/2).$

Now suppose $k\geq 2$
 and we have moved through Steps 2(b,c,d,e) for $k$. Hence, the algorithm has made it through these steps for $1\leq j \leq  k-1$ without halting. Thus: 
 \begin{enumerate}
\item[(b)] we have calculated $g^j(z')$ and $|Dg^j(z')|_e$, 
\item[(c)] concluded that the ball $B^e(g^j(z'),|Dg^j(z')|_e(3\sigma_N)(1+\gamma))$ is not disjoint from $\cN_2$, 
\item[(d)] concluded that the ball $B^e(g^j(z'),|Dg^j(z')|_e(4\sigma_N)(1-\gamma))$ does not entirely contain a reference ball from $\cN_2$, and
\item[(e)] showed  that $|Dg^j(z')|_e\sigma_N c(\gamma) \leq 2\beta'\sqrt{2}.$
 \end{enumerate}
Suppose  to the contrary that $g^{k-1}(z')\notin \cN(\cN_1,\delta/2)$.
We know $z'\in \cN(\cN_1,\delta/2)$, which is the base case $k=1$. Let $j\in \{1,\ldots,k-1\}$ be the smallest integer such that $z', g(z'), \ldots g^{j-1}(z') \in \cN(\cN_1,\delta/2)$ but $g^{j}(z')\notin \cN(\cN_1,\delta/2)$. 
Since $1\leq j\leq k-1$ we have the above statements (c), (d), and (e) are true for $j$. 

By (e) for $j$, 
$|Dg^j(z')|_e \sigma_N c(\gamma) \leq  2\beta'\sqrt{2},$
so $|Dg^j(z')|_e \sigma_N  \leq  2\beta'\sqrt{2}/(c(\gamma)).$
Combining this estimate with the distortion theorem, we have
$\frac{1}{2}{\rm diam}(g^j(\Omega_U)) \leq |Dg^j(z')|_e 4 \sigma_N (1+\gamma) 
\leq 4(1+\gamma) 2\beta'\sqrt{2}/(c(\gamma)). 
$
Using $\beta'<b$, 
Equation~\ref{eqn:defn-b} implies
$$
\frac{1}{2}{\rm diam}(g^j(\Omega_U)) \leq  8\sqrt{2}(1+\gamma) \beta'/(c(\gamma)) < \delta/2, 
$$
and thus $g^j(\Omega_U)$ is contained $\delta/2$-neighborhood of $g^j(z')$. If $g^{j}(z')\notin \cN(\cN_2,\delta/2),$ then $g^j(\Omega_U)$ is disjoint from $\cN_2$. But this means that the overestimate of $g^j(\Omega_U)$ is disjoint from $\cN_2$, so the overestimate of the subset $g^j(\Omega_L)$ must also be disjoint from $\cN_2$. This contradicts Step 2(c)'s conclusion for $j$. Hence, we have (iv): $g^j(z')\in \cN(\cN_2,\delta/2)$ for all $j\in \{1,\ldots, k-1\}$, so $g^{k-1}(z')\in \cN(\cN_2,\delta/2)$.  
\end{proof}

With these conditions established, we now show that the reporting of either $0$ or $1$ in the loop of Step 2 ((c),(d), or (e)) is correct. 

\begin{lemma} \label{lem:step2-correct}
     Let $p$ be a polynomial map of $\CC$ of degree $d\geq 2$. Beginning to apply Algorithm~\ref{alg:computeJ}, suppose we are examining $z'$ an ideal point of denominator $2^{N+2}$ for an $N\geq N'$.   Assume Algorithm~\ref{alg:computeJ} makes it to Step 2 but then halts for some $k\in \{1,\ldots,k_N\}$ in one of the Step 2(c), (d), or (e) tests, which leads to the reporting $0$ or $1$ for $z'$. Then this reporting is consistent with the definition of a $2^N$-approximation of $J$; i.e.,
    \begin{itemize}
     \item If \textbf{Step 2(c),(d), or (e)} reports $0$ for $z'$, then $d(z',J)\geq 2^{-N-2}$.
     \item If \textbf{Step 2(c),(d), or (e)} reports $1$ for $z'$, then $d(z',J) \leq 2^{-N-1}$.
    \end{itemize}
\end{lemma}

\begin{proof}
 Again, first note that if we make it to Step 2, then Step 1(a) did compute the neighborhoods $\cN_2, \cN_1$,  a consequence is that the map $p$ is hyperbolic, thus hyperbolicity is not required as an assumption.

Since we make it to Step 2, and $N\geq N'$ satisfying Equation~\ref{eqn:define_N'}, we know that (i), (ii), (iii), (iv) of Lemma~\ref{lem:step2-conditions} hold.
With these useful statements, we perform the loop for $k=1,\ldots,k_N$ through Steps 2(b, c, d, e), inductively and show each instance of reporting $0$ or $1$ is correct. 
For this proof, we apply the distortion theorem to $z'$ and to points in $\Omega_L$ that could be potential points in $J$ inside of $\cU'$. We use Proposition~\ref{prop:koebe-other-point} and Lemma~\ref{lem:proof-parts-1}. 

Stepping through the loop for $k\in \{1,\ldots, k_N\}$, we:

(b) calculate  $g^k(z')$ and $|Dg^k(z')|_e$. Then,

(c) test whether the ball $B^e(z',|Dg^k(z')|_e(3\sigma_N)(1+\gamma))$ is  disjoint from $\cN_2$. Since $z'\in \cN_1$, this ball is the overestimate of $g^k(\Omega_L)$. 

If it is, we reported $0$. Since we have here the overestimate of $g^k(\Omega_L)$ is disjoint from $\cN_2$, we know that $g^k(\Omega_L)$ is disjoint from $\cN_2$.
Since $J\subset \cN_2$ we know $g^k(\Omega_L)$ is disjoint from $J$, hence  $\Omega_L$ is disjoint from $J$. Thus, $d(z',J)> 2^{-N-2},$ the radius of $\Omega_L$, hence reporting $0$ is correct.

If we do not halt, 
we proceed to:

(d) test whether the ball $B^e(z',|Dg^k(z')|_e(4\sigma_N)(1-\gamma))$ contains a reference box from $\cN_2$. Since $z'\in \cN_1$, this ball is the underestimate of $g^k(\Omega_U)$. 

If it contains a reference box then we have reported $1$. Since every box of $\cN_2$ contains a point of $J$, since the underestimate of $g^k(\Omega_U)$ contains a reference box, this means $g^k(\Omega_U)$, hence $\Omega_U$, is guaranteed to contain a point of $J$. Thus $d(z',J)\leq 2^{-N-1}$ and reporting $1$ is permitted.

If not, we proceed to:

(e) test whether  $|Dg^k(z')|_e \cdot \sigma_N \cdot  c(\gamma) > 2\beta'\sqrt{2}.$

If this holds, we report $0$. We must show that is appropriate. 

Suppose to the contrary that there is a point $w'\in J \cap \cU'$. 
Then by (iii) of Lemma~\ref{lem:proof-parts-1}, since we found $|Dg^k(z')|_e\cdot \sigma_N \cdot c(\gamma) > 2\beta'\sqrt{2}$, we know  there is a reference box of $\cN_2$ inside the underestimate of $g^k(\Omega_U)$. 

  But this is a contradiction, as in Step 2(e) we are past 2(d), so no reference boxes are entirely contained in the underestimate of $g^k(\Omega_U)$. Thus there could not have been a point in $J$ in $\cU'$. 
   Hence, the reporting of $0$ for $z'$ is permissible.

On the other hand, if $|Dg^k(z')|_e \cdot \sigma_N \cdot  c(\gamma) \leq 2\beta'\sqrt{2},$ we pass Step 2(e) and are finished with the loop for this $k$. 
If $k<k_N$,  we go back up to Step 2(b) considering $k$ and advance to $k+1$. 
If $k=k_N$, and the algorithm did not halt, we have completed the loop and proceed to Step 2(f). Thus, the reporting of $0$ or $1$ for $z'$ in Steps 2(c),(d),(e) satisfies the conditions of a $2^{N}$-approximation. 
\end{proof}

Finally, we show that the algorithm is correct if it gets to the final step of the algorithm, Step 2(f), and then halts.

\begin{lemma}\label{lem:Step2f-correct}
          Let $p$ be a polynomial map of $\CC$ of degree $d\geq 2$. By applying Algorithm~\ref{alg:computeJ} we examine $z'$ an ideal point of denominator $2^{N+2}$ for an $N\geq N'$. Assume that  Algorithm~\ref{alg:computeJ} arrives at \textbf{Step 2(f)} and thus reports $0$. Then this is consistent with the definition of $2^N$-approximation of $J$, i.e.,
$d(z',J)\geq 2^{-N-2}$.
\end{lemma}

\begin{proof}
  Recall again that if we make it to Step 2, then Step 1(a) did halt and successfully computed the neighborhoods $\cN_2, \cN_1$, hence $p$ is hyperbolic so this doesn't need to be included in the assumptions. 
    
    If we make it to Step 2(f), none of the prior tests triggered a report. 
    
Now we need to show reporting $0$ in Step 2(f) is appropriate. Suppose to the contrary that we make it to \text{Step 2(f)}, but there is a point $w' \in J \cap \cU'$, where recall $\cU'=\overline{B(z',2^{-N-2})} \subset \Omega_L = B^e(z',3\sigma_N) $. 
Since $d_H(\Omega_U,\Omega_L)=\sigma_N$, we have $B^e(w',d(w',\partial \Omega_U ) ) \supset B^e(w',\sigma_N)=:W'$ independently of where $w'$ is located in $\cU'$.

As we loop through Steps 2(b, c, d, e),  we never report $0$ because $w' \in J\cap \cU'$ and by Lemma~\ref{lem:step2-correct}, the reporting of $0$ was valid in those steps. We would have reported $1$ for some $k$ in Step 2(d) if the underestimate for $g^k(\Omega_U)$ entirely contained a reference box. We now show  that $k_N$ was defined precisely so that this is true for some $k\in \{1,\ldots,k_N\}$ if there is a point in $J\cap \cU'$.

It follows from (iv) of Lemma~\ref{lem:step2-conditions} that if we are at Step 2(f), $g^j(z')\in \cN(\cN_2,\delta/2)$ for all $1\leq j\leq k_N$. Since $\cN_2\subset \cN_1$, we have 
$|Dg^{k_N}(z')|_e  \geq L^{k_N}$.
Hence Lemma~\ref{lem:proof-parts-1} implies:
$$
d_e(g^{k_N}(w'),g^{k_N}(w)) \geq d_e(w',w) \ c(\gamma) \  L^{k_N},
$$
for $w\in W'$, so $g^{k_N}(W') = g^{k_N}(B^e(w', \sigma_N)) \supset B^e ( g^{k_N}(w'), L^{k_N} \ c(\gamma) \ \sigma_N ).$

By definition of $k_N$ in Equation~\ref{eqn:defn-k-N} (also using that $\sigma_N = 2^{-N-3})$ we conclude that $k_N$ is the smallest positive integer that satisfies
$L^k   c(\gamma)  \sigma_N  > 2\beta'\sqrt{2} $. Putting these two inequalities together we have $\frac{1}{2}\text{diam}(g^{k_N}(W'))\geq \sigma_N c(\gamma) L^{k_N} > 2\beta'\sqrt{2},$ so since $|Dg^{k_N}(z')|_e  \geq L^{k_N}$ we have the hypothesis of (iii) of Lemma~\ref{lem:proof-parts-1} satisfied for $i=k_N$. Thus we know a reference box of $\cN_2$ is contained in the underestimate of $g^{k_N}(\Omega_U)$,  if there is a point $w'\in J\cap\cU'$. 

So  if there is a point $w'\in J\cap\cU'$, then for some $k\leq k_N$, there must be a reference box inside of the underestimate of $g^k(\Omega_U)$. But that would have triggered a report of $1$ in some step before 2(f). Since we got to Step 2(f), there could not have been a point $w'$ of $J$ in $\cU'$. Hence we may report $0$ at Step 2(f). 
\end{proof}

Combining the previous lemmas  establishes the following:
\begin{proposition}
\label{prop:correctness}
Let $p:\CC\to\CC$ be a hyperbolic polynomial of degree $d\geq 2$ and $N\in\bN$. Then Algorithm~\ref{alg:computeJ} produces a $2^{-N}$ approximation of $J$. 
\end{proposition}

The final two pieces of the proof of Theorem~\ref{thm:J-computable} are to show the algorithm halts, and to establish the running time. 

\begin{lemma}
    \label{lem:halting}Let $p:\CC\to\CC$ be a hyperbolic polynomial of degree $d\geq 2$.
    Then Algorithm~\ref{alg:computeJ} halts. 
\end{lemma}

\begin{proof}
Using the same arguments as in Lemma 3.9 of~\cite{BoydWolf-Skew1} we conclude that Step 1(a) of Algorithm~\ref{alg:computeJ} halts in finite time. Also,
we performed a finite, constant (independent of $N$ but dependent of $p$) number of calculations first, then for each $N$ and for a finite number of ideal points $z'$ with denominators $2^{N+2},$ we calculated the integer $k_N$ and only performed subsequent calculations up to at most the $k_N$-{th} iterate of the map or the derivative. Depending on the outcome of our calculations, we report $1$ or $0$, or if none of those conditions are met we arrive at Step 2(f) and  report $0$. This shows that the algorithm halts in finite time. 
\end{proof}

\begin{lemma}
    \label{lem:complexity}
    Let $p:\CC\to\CC$ be a hyperbolic polynomial of degree $d\geq 2$. Then Algorithm~\ref{alg:computeJ} produces a $2^N$-approximation of $J$ in time $O(N\cdot M(N))$, where $M(N)$ is the complexity of multiplying two $N$-bit digits. 
\end{lemma}

\begin{proof}
   The first computation in Step 1 produces a constant, independent of $N$, as it is carried out only once at the beginning and doesn't have to be repeated for higher $N$. 

Steps 1(d),(e) and (f) each require performing a test which takes a constant time independent of $N$. 

Step 2(a) is the calculation of $k_N$, which by Equation~\ref{eqn:defn-k-N} is a linear function of $N$. 
In Steps 2(b),(c),(d),(e) we perform a loop through $k$ from $1$ to $k_N$. 

In Step 2(b), we calculate up to the first $k_N$ points in the orbit of the ideal point $z'$ under $g$, and $|Dg^k(z')|$ for each $k$ from $1$ to at most $k_N$.
We note that, given the chain rule, some iterate derivatives may be stored to reduce the number of operations required, but this does not affect the asymptotic running time. 

In Step 2(c) we check if $B^e(z',|Dg^k(z')|\cdot 3\cdot 2^{-N-3}\cdot (1+\gamma))$ is disjoint from $\cN_2$. 

In Step 2(d), we perform a similar calculation, checking if $B^e(z',|Dg^k(z')|\cdot 4\cdot 2^{-N-3}\cdot (1-\gamma))$ contains any box of $\cN_2$. 

Step 2(e) is a simple check whether $|Dg^k(z')|_e \cdot 2^{-N-3}\cdot c(\gamma) > 2\beta'\sqrt{2}.$

The only aspects that affects the complexity in regards to $N$ is computing the points at higher precision for higher $N$. Thus, the complexity is a constant times $N$ times a constant multiple of the complexity of multiplying two $N$-bit digits which is $O(N \cdot M(N)).$
\end{proof}

Combining Proposition~\ref{prop:correctness} with Lemmas~\ref{lem:halting} and~\ref{lem:complexity} completes the proof of Theorem~\ref{thm:J-computable}.

\section{Concluding Remarks}
In this paper, we develop a new approach that re-establishes the poly-time computability of Julia sets for hyperbolic complex polynomials, and allows us to conclude lower semi-computability of the hyperbolicity locus for polynomials of fixed degree. The poly-time computability of hyperbolic Julia sets was originally established by Mark Braverman in 2005. Unlike Braverman’s result, which relies heavily on complex-analytic techniques such as the use of the Poincaré metric and Pick's theorem, our approach is purely dynamical in nature.

Specifically, we computationally establish uniform expansion in the Euclidean metric for a suitably identified iterate of the polynomial. Among other dynamical tools this step uses the characterization of the Julia set as the expanding chain recurrent set. We then apply several well-known dynamical tools, including invariance properties, quantitative shadowing, and distortion estimates. Since these dynamical techniques are readily available for conformal repellers, our results can be adapted to that setting. Moreover, it is our hope that the methods developed in this paper may serve as a foundation for future computational complexity results concerning more general dynamically defined invariant sets, such as hyperbolic saddle-sets, attractors and higher dimensional Julia sets.

\appendix 

\section{More Proof Details}

\begin{proof}[Proof of Proposition~\ref{prop:firststep}]
For completeness, we include here more details for this proof. 

First, given any positive integer $t$, we may consider the refined set $\cB'_{n+t}$: subdivide the boxes of $\cB'_n$ by placing a $(2^t)^2$-grid in each box of $\cB'_n$, then repeat the Algorithm 3.2, building edges where a box's (approximate but error-controlled) image overlaps another box, then refining to yield only boxes that lie in cycles and selecting only (but all of) the expanding components. Whenever we refine, since the $\cB_k$ sequence is nested, the smaller boxes are in one of the prior boxes so the expansion inequality doesn't change, so we have expansion by $1+\lambda'$ on $\cB'_{n+t}$. 

Now, due to this expansion, applying Corollary 2.19 of \cite{BoydWolf-Skew1} to $p^\nu$ on $\cB'_n$, we compute a radius $r'$ so that if two points in $\cB'_n$ are in a box of radius $\leq r'$, then the distance between their images under $p^\nu$ are pushed apart by at least $1+\lambda'$ times the distance between the two original points. By taking $r'$ slightly smaller, we may assume it's a dyadic rational. 
Now set 
\\ \centerline{
$
\beta := \min (r', \delta, (b-\ep_n)/2).
$}
\\This implies $\beta \leq r'< 2r'$ and $\ep_{n+t} + 2\beta \leq b$ for all $t\geq 0$, because $\ep_{n+t} \leq \ep_n$ (because the box diameter for boxes in a $\cB_k$ is $\ep_k = 2R/2^{k}$ where $R$ is determined by $p$). We will set $\beta' = \ep_{n+t}+2\beta$ for a $t$ to be calculated near the end of this proof, so we will have our $\cN_2$ be the union of boxes with centers in some $\cB'_{n+t}$ but with sidelength $\beta'$,  hence radii $\ep_{n+t}/2 + \beta$. Thus, we will have $\cN_2 \subset \cN(\cB'_n, \beta)$.
Now, the distance from the boundary of $\cB'_{n+t}$ to the boundary of $\cN_1$ is at least $2\delta$, because $d_H(\cB_n,\cN_1)=2\delta$ and $\cB_{n+t} \subset \cB_n$.
Since $\beta < \delta,$  we have the 
$\beta$-neighborhood of any box in $\cB'_{n+t}$ is contained in the $\delta$-neighborhood of $\cB'_n$, and hence at least a distance $\delta$ to $\cN_1$. Thus,
we have both 
$d_H(\cN_2,\cN_1)>\delta$ and 
$\cN(\cN_2, \delta) \subset \cN_1$, which will provide (iv). 




Since on $\cN_1$, $p^\nu$ is expanding by $L$, the same holds on subsets of $\cN_1$. Thus  for any $t\geq 0$, (iii) holds for $\cN(\cB'_{n+t}, \beta)$ hence will hold for our $\cN_2$.

Now we turn to finding the large enough $t$. We have $\cB'_k \subset \cR(\alpha_k)$, where $\alpha_k \downarrow 0$ as $k\to\infty$ and is computable. 
We use shadowing to find how close a real orbit is to our pseudo-orbits. 
Since $p^\nu$ is hyperbolic (because $p$ is), it has the shadowing property on each chain component.
We apply Equation (4.29) in \cite{Urbanski_book} to our setting. This provides a minimum (denoted by $\xi$ in \cite{Urbanski_book}) which in our setting is $2r'$ because the largest ball that fits inside of an image of a ball of radius $2r'$ has a radius bounded by $\onorm{Dp^\nu} 2r' 
 > 2r'$ since $\onorm{Dp^\nu} > 1$. 
By Proposition 4.3.4 in \cite{Urbanski_book}, we deduce that since $\beta < 2r'$ (by choice of $\beta$ above), we can compute a $t'$  sufficiently large that 
\\ \centerline{
$\alpha_{n+t'} < \min (2r', \lambda'\beta/4)).$}
\\Then if $t \geq t'$,   
there is a (unique) true orbit in $J$ that $\beta$-shadows any $\alpha_{n+t}$-pseudo-orbit. (For the reader who is carefully comparing this to \cite{BoydWolf-Skew1}, the $\alpha$ definition in the Type $J_2$-algorithm there included a $\lambda' \beta/2$ (though there $\beta$ was also $=2^N$). Here we have a $4$ in the denominator because we are only guaranteed expansion on our set of consideration by $\lambda'/2$, not $\lambda'$.)

Note this inequality also guarantees that $t'$ is so large that $\ep_{n+t} < \alpha_{n+t} \leq 2r'$, for any $t \geq t'$. We see below why this inequality is desired.  

If $t' > 0$, then we must refine $\cB'_n$ to yield a smaller boxes $\cB'_{n+t'}$, as described a bit above, and we'll consider $\cN(\cB'_{n+t'}, \beta)$. Again, box overlap is no problem.

To translate this back to box sizes, now with any $t\geq t'$, in any $\cB'_{n+t}$ with boxes of size $\ep_{n+t}$ which we know is $\leq 2r',$  we know by Corollary 3.5 of \cite{BoydWolf-Skew1} adapted to $p$ that the points in boxes of size $\ep_{n+t}$ are in pseudo-orbits of size $ \alpha_{n+t}$.
Thus $\beta$ defines the $t'$ sufficiently large that for $\alpha_{n+t} \leq  \min (2r', \lambda'\beta/4 )$, all $\alpha_{n+t}$-pseudo-orbits are $\beta$-shadowed by true orbits, and these are orbits in $J$ by choice of the chain component. 

Hence $\cB'_{n+t}$ lies within a $\beta$-neighborhood of $J$, i.e., $\cB'_{n+t} \subset \cN(J, \beta)$.
Thus, if (consistent with earlier notation) $\sB'_{n+t}$ is our collection of boxes whose union defines $\cB'_{n+t}$, then each box in $\cup_{B_j \in \sB'_{n+t}}\cN(B_j, \beta)$ contains a point of $J$, for all $t \geq t'$.

But we also know $J \subset \cB'_{n+t} \subset \cN(\cB'_{n+t}, \beta)$.
So, we have (i)
$
J \subset \cN(\cB'_{n+t}, \beta)$
and (ii) each box in $\cup \{ \cN(B_j, \beta): {B_j \in \sB'_{n+t}} \}$ contains a point of $J$ 
for any $t \geq t'$.

We established a few paragraphs above that (iii) holds for $\cN(\cB'_{n+t},\beta), \forall t \geq 0$.

It suffices to choose a $t\geq t'$ and set  $\cN_2 := \cup 
\{ \cN(B_j, \beta) : B_j \in \sB'_{n+t} \} =  \cN(\cB'_{n+t}, \beta)$, so $\beta':= \ep_{n+t}+2\beta.$ By definition of $\beta$ we then have $\beta' \in (0, b]$ as desired, (i),(ii),(iii) hold, and  $\cN(\cN_2,\delta) = \cN( \cB'_{n+t},\beta+\delta) \subset \cN_1$ as was explained a few paragraphs above, so (iv) holds. 

Note finally that $d_H(\cB'_{n+t},J) \leq {\ep_{n+t}}+\beta$ and $d_H(\cN_2,J) \leq \ep_{n+t}+2\beta=\beta'$.
\end{proof}

\bibliographystyle{amsplain}

\providecommand{\bysame}{\leavevmode\hbox to3em{\hrulefill}\thinspace}
\providecommand{\MR}{\relax\ifhmode\unskip\space\fi MR }
\providecommand{\MRhref}[2]{%
  \href{http://www.ams.org/mathscinet-getitem?mr=#1}{#2}
}

\end{document}